\documentclass[letter, oneside, 11pt, article]{memoir}

\usepackage[latin1]{inputenc}
\usepackage[T1]{fontenc}
\usepackage{ae}
\usepackage{aecompl}

\usepackage{amssymb}
\usepackage{amsmath}
\usepackage{amsfonts}
\usepackage{amsthm}

\usepackage[english]{babel}

\usepackage{graphicx}

\usepackage[plain, english]{fancyref}

\usepackage{bbm}
\usepackage{mathrsfs}

\usepackage{stmaryrd}

\usepackage{ifthen}

\usepackage{xspace}

\usepackage{calc}
\usepackage{color}

\usepackage[ps,all,dvips]{xy}
\SelectTips{cm}{12}
\CompileMatrices

\newcommand{\mathset}[1]{\mathbbm{#1}}

\newcommand{\setZ}{\mathset{Z}}

\newcommand{\setR}{\mathset{R}}
\newcommand{\setC}{\mathset{C}}

\renewcommand{\d}{\partial}

\renewcommand{\phi}{\varphi}
\renewcommand{\epsilon}{\varepsilon}
\renewcommand{\otimes}{\varotimes}

\newcommand{\sotimes}{\mathbin{\raise1.5pt\hbox{
      $\scriptscriptstyle\otimes$}}}

\newcommand{\boldnabla}{\mbox{\boldmath$\nabla$}}

\DeclareMathOperator{\im}{Im}       
\DeclareMathOperator{\id}{Id}       
\DeclareMathOperator{\Tr}{Tr}
\DeclareMathOperator{\End}{End}
\DeclareMathOperator{\Id}{Id}

\theoremstyle{plain}
\newtheorem{lemma}{Lemma}
\newtheorem{proposition}[lemma]{Proposition}
\newtheorem{theorem}[lemma]{Theorem}
\newtheorem{corollary}[lemma]{Corollary}
\newtheorem{definition}[lemma]{Definition}

\theoremstyle{definition}

\makeatletter       

\let\@xp\expandafter 
\newcommand\DefineFancyrefPrefix[2]{%
  \@namedef{fancyref#1labelprefix}{#1}%
  \@namedef{Fref#1name}{#2}%
  \@namedef{fref#1name}{\MakeLowerCase{\@nameuse{Fref#1name}}}%
  \def\@style{vario}%
  \@xp\@xp\@xp\frefformat\@xp\@xp\@xp\@style\@xp\csname
  fancyref#1labelprefix\endcsname
  {%
    \@nameuse{fref#1name}\fancyrefdefaultspacing##1##3%
  }%
  \@xp\@xp\@xp\Frefformat\@xp\@xp\@xp\@style\@xp\csname
  fancyref#1labelprefix\endcsname
  {%
    \@nameuse{Fref#1name}\fancyrefdefaultspacing##1##3%
  }%
  \def\@style{plain}%
  \@xp\@xp\@xp\frefformat\@xp\@xp\@xp\@style\@xp\csname
  fancyref#1labelprefix\endcsname
  {%
    \@nameuse{fref#1name}\fancyrefdefaultspacing##1%
  }%
  \@xp\@xp\@xp\Frefformat\@xp\@xp\@xp\@style\@xp\csname
  fancyref#1labelprefix\endcsname
  {%
    \@nameuse{Fref#1name}\fancyrefdefaultspacing##1%
  }%
}
\makeatother

\DefineFancyrefPrefix{ax}{Axiom}
\DefineFancyrefPrefix{cor}{Corollary}
\DefineFancyrefPrefix{lem}{Lemma}
\DefineFancyrefPrefix{prop}{Proposition}
\DefineFancyrefPrefix{thm}{Theorem}
\DefineFancyrefPrefix{def}{Definition}
\DefineFancyrefPrefix{rem}{Remark}

\makepagestyle{niels}

\makeoddhead{niels}{}{}{\thepage}

\chapterstyle{article}
\begin{document}

\title{Hitchin's Connection in Half-form Quantization} \author{Jørgen
  Ellegaard Andersen, Niels Leth Gammelgaard \\ and \\Magnus Roed
  Lauridsen}

\maketitle

\begin{abstract}
  We give a differential geometric construction of a connection in the
  bundle of quantum Hilbert spaces arising from half-form corrected
  geometric quantization of a prequantizable, symplectic manifold,
  endowed with a rigid, family of K\"ahler structures, all of which
  give vanishing first Dolbeault cohomology groups.

  In \cite{A1} Andersen gave an explicit construction of
  Hitchin's connection in the non-corrected case using additional
  assumptions. Under the same assumptions we also give an explicit
  solution in terms of Ricci potentials.  Morover we show that if
  these are carefully chosen the construction coincides with the
  construction of Andersen in the non-corrected case.
\end{abstract}

\chapter{Introduction}

Hitchin constructed in \cite{MR1065677} a connection over Teichmüller
space. This Hitchin connection is a connection in the bundle obtained
from geometric quantization of the moduli spaces of flat
$\mathrm{SU(n)}$-connections on a closed oriented surface.  The significance of
this connection is its relation to $(2+1)$-dimensional
Reshetikhin-Turaev TQFT (\cite{MR1036112} and \cite{MR1091619}). In
fact, this geometric construction of these TQFT's was proposed by
Witten in \cite{MR990772}, where he derived, via the Hamiltonian
approach to quantum Chern-Simons theory, that the geometric
quantization of the moduli spaces of flat connections should give the
two dimensional part of the theory. Further, he proposed an
alternative construction of the two dimensional part of the theory via
WZW-conformal field theory. This theory has been studied
intensively. In particular the work of Tsuchiya, Ueno and Yamada in
\cite{MR1048605} provided the major geometric constructions and
results needed. In \cite{MR1797619}, their results was used to show
that the category of integrable highest weight modules of level $k$
for the affine Lie algebra associated to any simple Lie algebra is a
modular tensor category. Further in \cite{MR1797619} this result is
combined with the work of Kazhdan and Lusztig (\cite{MR1186962},
\cite{MR1239506} and \cite{MR1239507}) and the work of Finkelberg
\cite{MR1384612} to argue that this category is isomorphic to the
modular tensor category associated to the corresponding quantum group,
from which Reshetikhin and Turaev constructed their
TQFT. Unfortunately, these results do not allow one to conclude the
validity of the geometric constructions of the two dimensional part of
the TQFT proposed by Witten.  However, in joint work with Ueno,
\cite{AU4}, \cite{AU3}, \cite{MR2339577} and \cite{MR2306213}, we have
given a proof, based mainly on the results of \cite{MR1048605}, that
the TUY-construction of the WZW-conformal field theory after twist by
a fractional power of an abelian theory, satisfies all the axioms of a
modular functor.  Furthermore, we have proved that the full
$2+1$-dimensional TQFT that results from this is isomorphic to the one
constructed by BHMV via skein theory mentioned above. Combining this
with the Theorem of Laszlo \cite{MR1669720}, which identifies
(projectively) the representations of the mapping class groups one
obtains from the geometric quantization of the moduli space of flat
connections with the ones obtained from the TUY-constructions, one
gets a proof of the validity of the construction proposed by Witten in
\cite{MR990772}.

In \cite{MR1100212}, Axelrod, Della Pietra and Witten gave a
differential geometric construction of the Hitchin connection by using
a method of symplectic reduction from the infinite dimensional space
of all $\mathrm{SU(n)}$-connections. In \cite{A1} Andersen constructed the
Hitchin connection in a more general setting. A corollary of the
results in \cite{A1} is that the connection constructed by Axelrod,
Della Pietra and Witten in \cite{MR1100212} is the same as Hitchin's
connection constructed in \cite{MR1065677}.

In this paper, we extend the setting from \cite{A1}, in which we can
construct the Hitchin connection. Let us describe this setting.

Let $(M,\omega)$ be a symplectic manifold of dimension $2m$, and let
$J$ be a family of K\"ahler structures on $M$, parametrized smoothly
by a manifold $\mathcal{T}$. Along any vector field $V$ on
$\mathcal{T}$, we can differentiate $J$ to get a map $V[J]\colon
\mathcal{T} \to C^\infty(M, \End(TM_{\setC}))$.

Define $\tilde{G}(V)\in C^\infty(M,S^2 (TM_{\setC}))$ by
\begin{align*}
  V[J]=\tilde{G}(V)\omega.
\end{align*}
Letting $T_\sigma$ denote the holomorphic tangent bundle on $(M,
J_\sigma)$ for any $\sigma \in \mathcal{T}$, we can further define
$G(V) \in C^\infty(M,S^2 (T))$ by the equation
\begin{align*}
  \tilde{G}(V)=G(V)+\bar G(V),
\end{align*}
for all real vector field $V$ on $\mathcal{T}$. We shall assume that
the family $J$ is \emph{rigid}, in the sense of \Fref{def:4}, meaning
that $G(V)_\sigma$ is a holomorphic section of $S^2 (T_\sigma)$. 

In case the second Stiefel-Whitney class vanishes, we can choose a
metaplectic structure on $(M, \omega)$, which gives rise to a choice
of a square root $\delta_\sigma$ of the canonical line bundle
$K_\sigma \to M_\sigma$, varying smoothly in the parameter $\sigma \in
\mathcal{T}$.

Now assume that $(M,\omega)$ is prequantizable in the sense that there
exists a Hermitian line bundle $\mathcal{L}$ with a compatible
connection $\nabla^\mathcal{L}$ of curvature
\begin{align*}
  R_{\nabla^\mathcal{L}} = -i \omega.
\end{align*}
The Levi-Civita connection $\tilde \nabla_\sigma$, corresponding to
the K\"ahler metric on $M_\sigma$, induces a connection in the line
bundle $\delta_\sigma\to M_\sigma$, and thus we get a connection
$\nabla_\sigma$ in $L^k\otimes \delta_\sigma\to M_\sigma$ giving this
bundle the structure of a holomorphic line bundle.

For every $\sigma \in \mathcal{T}$, we have the infinite dimensional
vector space $\mathcal{H}_\sigma^{(k)} = C^\infty(M, \mathcal{L}^k
\otimes \delta_\sigma)$, and we consider subspace of holomorphic
sections
\begin{equation*}
  H_\sigma^{(k)} = H^0(M_\sigma, \mathcal{L}^k\otimes\delta_\sigma) = \{ s\in
  C^\infty(M, \mathcal{L}^k\otimes\delta_\sigma) \mid
  \nabla_\sigma^{0,1}s = 0 \}.
\end{equation*}
It is not clear that the spaces $\mathcal{H}^{(k)}_\sigma$ form a
smooth vector bundle $\mathcal{H}^{(k)} \to \mathcal{T}$. However, it
is a corollary of our construction that, under the assumptions stated
in \Fref{thm:1}, the spaces $\mathcal{H}^{(k)}_\sigma$ indeed form a
smooth bundle over $\mathcal{T}$ and that $H^{(k)} \to \mathcal{T}$ is
a smooth subbundle of $\mathcal{H}^{(k)}$.

In the bundle $C^\infty(M, T) \to \mathcal{T}$, we have a connection
$\hat \nabla^T$ defined by the formula
\begin{align*}
  \hat \nabla^T_V \zeta = \pi^{1,0} V[\zeta],
\end{align*}
where $\pi^{1,0}_\sigma \colon TM_\setC \to T_\sigma$ is the
projection, and $V[\zeta]$ denotes differentiation in the trivial
bundle $\mathcal{T} \times C^\infty(M, TM_\setC)$ (see section
\ref{RC} for further details). This induces a connection in
$C^\infty(M, \delta) \to \mathcal{T}$, and with the help of the
trivial connection in $\mathcal{T} \times C^\infty(M, \mathcal{L}^k)$
this induces a connection $\hat \nabla^r$ in $\mathcal{H}^{(k)} \to
\mathcal{T}$, called the reference connection.

\begin{theorem}
  \label{thm:1}
  Let $(M, \omega)$ be a prequantizable, symplectic manifold with
  vanishing second Stiefel-Whitney class. Further let $J$ be a rigid
  family of K\"ahler structures on $M$ all satifying $H^{0,1}(M) =
  0$. Then for any vector field $V$ on $\mathcal{T}$, the Hitchin
  connection $\boldnabla$ in the bundle $H^{(k)}$ is given by
  \begin{align*}
    \boldnabla_V = \hat \nabla^r_V + \frac{1}{4k}(\Delta_{G(V)} +
    H(V)),
  \end{align*}
  where $\hat \nabla^r$ is the reference connection, $\Delta_{G(V)}$
  is the second order differential operator $\Delta_{G(V)} = \Tr
  \nabla G(V) \nabla$, and $H$ is any one-form on $\mathcal{T}$, with
  values in $C^\infty(M)$, satisfying $\bar \d_M H(V) = \frac{i}{2}
  \Tr \tilde \nabla (G(V) \rho)$. Such a one-form $H$ exists and is
  unique up to addition of the pullback of an ordinary one-form on
  $\mathcal{T}$.
\end{theorem}

In fact, we can consider rigidity as a condition on the vector fields
$V$, rather than considering it as a condition on the family of
K\"ahler structures.  We will then get a partial connection, which is
defined in these rigid directions. This shows that as soon as the
family $\mathcal{T}$ contains points $\sigma$ such that
$H^0(M_\sigma, S^2(T)) \neq 0$, then the constructions provide a partial
connection defined at least on some non-zero tangential directions in
$\mathcal{T}$.

Assume now that $M$ is compact with $H^1(M, \setR) = 0$. Notice, that
by the Hodge decomposition theorem, $H^{0,1}(M) = 0$ for any K\"ahler
structure on $M$.

Now we further assume, that $\mathcal{T}$ is a complex manifold, and
that $J$ is a \emph{holomorphic} family in the sense of \Fref{def:3},
or equivalently that $J$ gives rise to a complex structure on
$\mathcal{T} \times M$.

We then consider the non-corrected setting of geometric quantization
of $(M, \omega)$, namely
\begin{equation*}
  \tilde H_\sigma^{(k)} = H^0(M_\sigma, \mathcal{L}^k) = \{ s \in
  C^\infty(M, \mathcal{L}^k) \mid
  (\nabla^\mathcal{L})_\sigma^{0,1}s = 0 \}.
\end{equation*}
Under the additional assumption, that the real first Chern class of $(M,
\omega)$ is given by
\begin{align}
  \label{eq:1}
  c_1(M, \omega) = n[\tfrac{\omega}{2\pi}], \qquad n \in \setZ,
\end{align}
there is a construction, due to Andersen (\cite{A1}), of a Hitchin
connection in the trivial bundle $\mathcal{T}\times C^\infty(M,
\mathcal{L}^k)$ over $\mathcal{T}$, which preserves the subbundle
$\tilde H^{(k)} \to \mathcal{T}$, extending Hitchin's connection
constructed in \cite{MR1065677}. Now when \eqref{eq:1} is satisfied,
we are able to give an explicit formula for the one-form $H$, namely
\begin{align*}
  H(V) = - 2nV'[F] - \d_M F G(V) \d_M F - \Tr \tilde \nabla (G(V) \d_M
  F),
\end{align*}
where $F$ is any smooth family of Ricci potentials on $M$ over
$\mathcal{T}$. Moreover the following theorem says, that if we choose
the right normalization of the Ricci potentials, we can compare the
Hitchin connection given by \Fref{thm:1} with the one constructed by
Andersen in the non-corrected case, and in fact they agree.

\begin{theorem}
  \label{thm:2}

  Let $(M, \omega)$ be a compact, prequantizable, symplectic manifold
  with vanishing second Stiefel-Whitney class, and
  $H^1(M,\setR)=0$. Further let $J$ be a rigid, holomorphic family of
  K\"ahler structures on $M$ parametrized by a complex manifold
  $\mathcal{T}$. Assume that the first Chern class of $(M, \omega)$
  is divisible by an integer $n$ and that its image in $H^2(M,\setR)$ satisfies
  \begin{align*}
    c_1(M, \omega) = n[\tfrac{\omega}{2\pi}].
  \end{align*}
  Then around every point $\sigma \in \mathcal{T}$, there exists an
  open neighbourhood $U$, a local smooth family $\tilde F$ of Ricci
  potentials on $M$ over $U$ and an isomorphism of vector bundles over
  $U$
  \begin{align*}
    \phi \colon \tilde H^{(k-n/2)} \vert_U \to H^{(k)} \vert_U,
  \end{align*}
  such that
  \begin{align*}
    \phi^* \boldnabla = \tilde \boldnabla,
  \end{align*}
  where $\phi^* \boldnabla$ is the pullback of the Hitchin connection
  given by \Fref{thm:1}, and $\tilde \boldnabla$ is the Hitchin
  connection in $\tilde H^{(k-n/2)}$ constructed in \cite{A1}, both of
  which are expressed in terms of $\tilde F$.
\end{theorem}

We plan to address the computation of the curvature and removal of the
rigidity condition in a forthcoming publication.  Also, we find it
interesting to analyze the relation between the connection constructed
in this paper and the "$L^2$-induced" constructed by Charles in
\cite{MR2276452}.  Futher we intend to consider this new construction
in the moduli space setting, in which Hitchin originally constructed
his connection, and which was applied further by Andersen in
\cite{MR2195137}.

Further, we find it very interesting to explore the role of Toeplitz
operators and their relation to the Hitchin connection constructed in
the general setting considered in this paper. In particular it would
be interesting to under stand if the results in \cite{A1}, \cite{A2}
and \cite{A3} can be generalized to this setting. For the first steps
in this direction see also \cite{A4}.

This paper is organized as follows. In section 2 we introduce
half-form corrected geometric quantization and the notion of
metaplectic structure. Section 3 is devoted the reference connection
and the calculation of its curvature. In section 4 we derive an
equation that the Hitchin connection should satisfy. Then we give a
solution to this equation and prove \Fref{thm:1}. Finally, in section
5, we study the relation between our construction and the construction
of \cite{A1} in the non-corrected case, culminating with a proof of
\Fref{thm:2}.

\chapter{Half-form Quantization and Metaplectic Structure}

Consider an almost complex structure $J$ on $M$, which is
compatible with the symplectic structure in the sense that
\begin{align*}
  g_J(X,Y) = \omega(X, JY)
\end{align*}
defines a Riemannian metric on $M$. We shall denote the resulting
Riemannian manifold by $M_J$.

The almost complex structure $J$ induces a splitting
\begin{align*}
  TM_\setC = T_J \oplus \bar T_J
\end{align*}
of the complexified tangent bundle into the eigenspaces of $J$
corresponding to the eigenvalues $i$ and $-i$ respectively. This
splitting is explicitly given by the projections onto each summand
\begin{align}
  \label{eq:2}
  \begin{aligned}
    \pi^{1,0}_J &= \tfrac{1}{2}(\Id - iJ) \qquad
    T_J = \im (\pi^{1,0}_J) \\
    \pi^{0,1}_J &= \tfrac{1}{2}(\Id + iJ) \qquad \bar T_J = \im
    (\pi^{0,1}_J).
  \end{aligned}
\end{align}
The fact that $T_J$ and $\bar T_J$ are the eigenspaces of $J$,
corresponding to the eigenvalues $i$ respectively $-i$, is easily
verified from these formulas. Very often we shall use the notation $X'
= \pi^{1,0}_J X$ and $X'' = \pi^{0,1}_J X$ for vector fields $X$ on
$M$.

Tensors, such as the symplectic form and associated metric, are
extended complex linearly to $TM_\setC$.

We recall that the first
Chern class $c_1(M_J)$ is equal to minus the first Chern class of the
canonical line bundle
\begin{align*}
  K_J = \bigwedge^m T_J^*.
\end{align*}
By integrality, $c_1(M_J)$ is independent of $J$ since the space of
compatible almost complex structures on $(M,\omega)$ is
contractible. Thus, the first Chern class is an invariant of the
symplectic manifold rather than the almost complex one.

Let us assume, that the second Stiefel-Whitney class $w_2(M)$
vanishes. Since the reduction modulo 2 of the first Chern class, that
is the image of $c_1(M)$ under the map $H^2(M, \setZ) \to H^2(M,
\setZ_2)$, is equal to the second Stiefel-Whitney class, this implies
that the first Chern class of $M$ is even. Thus the fist Chern class
of $K_J$ is even, which is equivalent to the existence of a square
root $\delta_J$ of $K_J$. We shall see later that the choice of such a
$\delta_J$ determines a square root of the canonical line bundle for
every other almost complex structure on $M$.

The metric on $M_J$ gives rise to the Levi-Civita connection $\tilde
\nabla_J$. As usual we get an induced metric and compatible connection
in all tensor bundles over $M$, and we shall denote all of these by
$g_J$ and $\tilde \nabla_J$ as well.

The metric also induces a Hermitian structure $h^T_J$ in $T_J$ given
by
\begin{align*}
  h^T_J(X, Y) = g_J(X, \bar Y),
\end{align*}
for any vectors $X$ and $Y$ in $T_J$.  If we further assume that $J$
is parallel, with respect to the Levi-Civita connection $\tilde
\nabla_J$, then $J$ must be integrable and $M_J$ K\"ahler. In this
case $\tilde \nabla_J$ preserves the holomorphic tangent bundle $T_J$
inducing a connection $\nabla^T_J$ compatible with $h^T_J$. These in
turn induce a Hermitian structure $h^K_J$ and compatible connection
$\nabla^K_J$ in the canonical line bundle $K_J$.

The Ricci tensor $r_J$ on $M_J$ is given by the following trace of the
K\"ahler curvature
\begin{align*}
  r_J(X,Y) = \Tr ( Z \mapsto \tilde R(Z, X)Y ),
\end{align*}
and the Ricci form $\rho_J$ is the associated (1,1)-form given by
\begin{align*}
  \rho_J(X,Y) = r(JX, Y).
\end{align*}
We recall for future use that the canonical line bundle $K_J$ has
curvature $i\rho_J$.

Finally $h^K_J$ and $\nabla^K_J$ induce a Hermitian strucuture
$h^\delta_J$ and compatible connection $\nabla^\delta_J$ in the line
bundle $\delta_J$.

\begin{definition}
  \label{def:1}
  A \emph{prequantum line bundle} over the symplectic manifold
  $(M,\omega)$ is a Hermitian line bundle $\mathcal{L}$ with a
  compatible connection $\nabla^\mathcal{L}$ of curvature
  \begin{align*}
    R_{\nabla^\mathcal{L}} = -i \omega,
  \end{align*}
  where $R_\nabla(X,Y) = [\nabla_X, \nabla_Y] - \nabla_{[X,Y]}$. Such
  a triple $(\mathcal{L}, h^\mathcal{L}, \nabla^\mathcal{L})$ is
  denoted a prequantum line bundle, and we say that the symplectic
  manifold is \emph{prequantizable} if it admits such a bundle.
\end{definition}

Evidently, a necessary condition for the existence of a prequantum
line bundle is that the class $[\tfrac{\omega}{2\pi}]$ in $H^2(M,
\setR)$ is integral, and in fact this is also sufficient. Moreover,
inequivalent choices of prequantum line bundles are parametrized by
the first cohomology $H^1(M, U(1))$ with coefficients in the circle
group $U(1) \subset \setC$ (see for instance \cite{MR1183739}). We
shall assume that $M$ is prequantizable, and fix a prequantum line
bundle $(\mathcal{L}, h, \nabla^\mathcal{L})$.

Now $h^\mathcal{L}$ and $h^\delta_J$ induce a Hermitian structure
$h_J$ in the line bundle $\mathcal{L}^k \otimes \delta_J$, and we have
a compatible connection $\nabla_J$, induced by $\nabla^\mathcal{L}$
and $\nabla^\delta_J$. Since $\mathcal{L}^k \otimes \delta_J$ has
curvature $-ik\omega + \frac{i}{2} \rho_J$, which is of type (1,1), the
operator
\begin{align*}
  \nabla_J^{0,1} = \pi^{0,1}_J \nabla_J
\end{align*}
defines a $\bar \d$-operator in $\mathcal{L}^k \otimes \delta_J$,
making this a holomorphic line bundle over $M_J$ (see
e.g. \cite{MR702806}). If we consider the space $\mathcal{H}^{(k)}_J =
C^\infty(M, \mathcal{L}^k \otimes \delta_J$) of smooth sections, then
the operator $\nabla_J^{0,1}$ gives rise to the subspace $H^{(k)}_J$
of holomorphic sections
\begin{align*}
  H^{(k)}_J = H^0(M_J, \mathcal{L}^k \otimes \delta_J) = \{ s \in
  C^\infty(M_J,\mathcal{L}^k \otimes \delta_J) \mid \nabla^{0,1}_J s =
  0 \}.
\end{align*}
We can define a Hermitian inner product on this space by
\begin{align*}
  \langle s_1 , s_2\rangle = \frac{1}{m!} \int_M h^\mathcal{L}(s_1,
  s_2) \omega^m,
\end{align*}
and if we consider the space of square integrable functions we obtain
a Hilbert space. This is the Hilbert space resulting from the
half-form corrected geometric quantization of the K\"ahler manifold
$M_J$.\\

We will construct a connection in $\mathcal{H}^{(k)}$ and prove, that
under certain conditions this connection preserves the infinitesimal
condition for being contained in the subspaces $H^{(k)}_J$. From this
we conclude, that the spaces $H^{(k)}_J$ form a vector bundle over a
manifold that parametrizes choices of $J$, and the fibers $H^{(k)}_J$
are related using parallel translation of the induced connection,
which we will call the Hitchin connection.

To be able to do this, we should pay closer attention to the way we
choose the half-form bundle $\delta_J$. Clearly, there is more than
one choice of a square root of $K_J$ (when it exists), and we would
like to choose $\delta_J$ in a unified way for different $J$. This is
were the notion of a metaplectic structure comes into the picture.

Consider the positive Lagrangian Grassmannian $L^+M$ consisting of
pairs $(p, J_p)$, where $p \in M$ and $J_p$ is a compatible almost
complex structure on the tangent space $T_pM$. This space has the
structure of a smooth bundle over $M$, with the obvious projection,
and with sections corresponding precisely to almost complex structures
on $M$.

At each point $(p, J_p) \in L^+M$, we can consider the one dimensional
space $K_{J_p} = \bigwedge^m T^*_{J_p}$. These form a smooth bundle
$K$ over $L^+M$, and the pullback by a section of $L^+M$ yields the
canonical line bundle associated to the almost complex structure on
$M$ given by the section.

We want to find a square root $\delta \to L^+M$ of the bundle $K \to
L^+M$. Such a square root is called a metaplectic structure on
$M$. Since $L^+M$ has contractible fibers, we can find local
trivializations of $K$ with constant transition functions along the
fibers. The construction of a metaplectic structure on $M$ amounts to
choosing square roots of these transition functions in such a way that
they still satisfy the cocycle conditions. But since the transition
functions are constant along the fibers, we only have to choose a
square root at a single point in each fiber. In other words, a square
root $\delta_J$ of $K_J$, for a single almost complex structure $J$ on
$M$, determines a metaplectic structure. We summarize this in a
proposition.

\begin{proposition}
  \label{prop:1}
  Let $M$ be a manifold with vanishing second Stiefel-Whitney class,
  and let $\omega$ be any symplectic structure on $M$. Then $(M,
  \omega)$ admits a metaplectic structure $\delta \to L^+M$.
\end{proposition}

For the rest of this paper, we shall assume that $M$ satisifies the
conditions of this proposition, and fix a metaplectic structure
$\delta$. In this way, for every almost complex structure $J$ on $M$,
viewed as a section of $L^+M$, we have a canonical choice of square
root of the canonical line bundle, given as the pullback of $\delta$
by $J$. \\

\chapter{The Reference Connection}\label{RC}

Returning to the setup of the introduction, consider a manifold
$\mathcal{T}$, and assume that we have a smooth family $J \colon
\mathcal{T} \to C^\infty(M, \End(TM))$ of K\"ahler structures on $M$,
parametrized by $\mathcal{T}$. More precisely $J$ is a smooth section
of the pullback bundle $\pi_M^* \End(TM) \to \mathcal{T} \times M$,
where $\pi_M \colon \mathcal{T} \times M \to M$ is the projection,
such that for every $\sigma \in \mathcal{T}$, the endomorphism
$J_\sigma$ defines a complex structure on $M$, turning this into a
K\"ahler manifold $M_\sigma$. As in the previous section, the K\"ahler
metric is given by
\begin{align*}
  g_\sigma(X,Y) = \omega(X,J_\sigma Y),
\end{align*}
and $J_\sigma$ induces a splitting $TM_\setC = T_\sigma \oplus \bar
T_\sigma$. Also we write $X'_\sigma = \pi^{1,0}_\sigma X$ and
$X''_\sigma = \pi^{0,1}_\sigma X$ for any vector field $X$ on $M$.

Viewing the family $J$ as a map $\mathcal{T} \times M \to L^+M$, we
get a smooth bundle $\delta \to \mathcal{T} \times M$, by pulling back
the metaplectic structure on $M$.  For any $\sigma \in \mathcal{T}$,
the restriction
\begin{align*}
  \delta_\sigma = \delta \vert _{\{\sigma\} \times M} \to M
\end{align*}
is a square root of the canonical line bundle $K_\sigma$ on
$M_\sigma$. Moreover the Hermitian structure $h^\delta_\sigma =
h^\delta_{J_\sigma}$ in $\delta_\sigma$ gives rise to a Hermitian
structure $h^\delta$ on $\delta$.
Let $\pi_M : \mathcal{T} \times M \to M$ denote the projection and
define
\begin{align*}
  \mathcal{\hat L} = \pi_M^* \mathcal{L} = \mathcal{T} \times
  \mathcal{L},
\end{align*}
with Hermitian metric $\hat h^\mathcal{L} = \pi_M^*
h^\mathcal{L}$. When objects are extended to the product $\mathcal{T}
\times M$, we shall often use a hat to indicate, that we are dealing
with the extended object. Then $\mathcal{\hat L} \otimes \delta$
becomes a smooth line bundle over $\mathcal{T} \times M$ with
Hermitian metric $\hat h$ induced by $\hat h^\mathcal{L}$ and
$h^\delta$.

As in the previous section we consider the space
$\mathcal{H}^{(k)}_\sigma = C^\infty(M_\sigma, \mathcal{L}^k \otimes
\delta_\sigma)$, in which the connection $\nabla_{J_\sigma}$, which we
shall denote by $\nabla_\sigma$, gives rise to the subspace of
holomorphic sections
\begin{align*}
  H^{(k)}_\sigma= H^0(M_\sigma, \mathcal{L}^k \otimes \delta_\sigma) =
  \{ s \in \mathcal{H}^{(k)}_\sigma \mid \nabla^{0,1}_\sigma s = 0 \}.
\end{align*}
In fact the spaces $\mathcal{H}^{(k)}_\sigma$ form a smooth vector
bundle $\mathcal{H}^{(k)}$ over $\mathcal{T}$. We will
construct a connection in $\mathcal{H}^{(k)}$ which preserves the
spaces $H^{(k)}_\sigma$, thereby proving that these form a smooth
 subbundle $H^{(k)}$ of $\mathcal{H}^{(k)}$, and at the same
time giving a connection in $H^{(k)}$. \\

First we define a connection $\hat \nabla^\mathcal{L}$ in
$\mathcal{\hat L}$ simply by extending $\nabla^\mathcal{L}$ as the
trivial connection in directions tangent to $\mathcal{T}$, i.e.  $\hat
\nabla^\mathcal{L}$ is the pullback connection in the pullback
bundle $\mathcal{\hat L}$. Concretely, if $X$ is a vector field on
$\mathcal{T} \times M$, which is tangent to $M$, and $s$ is a section
of $\mathcal{\hat L}$, then we define
\begin{align*}
  (\hat \nabla^\mathcal{L}_X s)_{(\sigma, p)} = (\nabla^\mathcal{L}_X
  s_\sigma)_p.
\end{align*}
For any vector field $V$ on $\mathcal{T} \times M$, which is tangent
to $\mathcal{T}$, we have that
\begin{align*}
  (\hat \nabla^\mathcal{L}_V s)_{(\sigma, p)} = V[s_p]_\sigma.
\end{align*}
Here $V[s_p]_\sigma$ denotes differentiation at $\sigma \in \mathcal{T}$ along $V$ of
$s_p,$ as a section of the trivial bundle $\mathcal{T} \times
\mathcal{L}_p$.

Now $\hat \nabla^\mathcal{L}$ is easily seen to be compatible with the
Hermitian structure $\hat h^\mathcal{L}$, and for future reference we
give the curvature, which is easily calculated.
\begin{lemma}
  \label{lem:1}
  The curvature of $\hat \nabla^\mathcal{L}$ is given by
  \begin{align*}
    R_{\hat \nabla^\mathcal{L}} = \pi_M^* R_{\nabla^\mathcal{L}} = -i\pi_M^* \omega,
  \end{align*}
  where $\pi_M \colon \mathcal{T} \times M \to M$ denotes the
  projection.
\end{lemma}

Next we define a connection $\hat \nabla^T$ in the bundle $T \to
\mathcal{T} \times M$ in the following way. In the directions tangent
to $M$, simply take $\hat \nabla^T$ to be the connection $\nabla^T$
induced from the Levi-Civita connection. More explicitly we define,
for any section section $Y$ of $T$ and any vector $X \in T_pM$,
\begin{align}
  \label{eq:3}
  (\hat \nabla^{T}_X Y)_{(\sigma,p)} = ((\nabla^T_\sigma)_X
  Y_\sigma)_p,
\end{align}
where $\nabla^T_\sigma$ denotes $\nabla^T_{J_\sigma}$. For the
directions along $\mathcal{T}$, we let $V \in T_\sigma \mathcal{T}$ be
any vector on $\mathcal{T}$ and define
\begin{align}
  \label{eq:4}
  (\hat \nabla^T_V Y)_{(\sigma,p)} = \pi_\sigma^{1,0} V[Y_p]_\sigma,
\end{align}
for any section $Y$ of $T$, where $V[Y_p]$ denotes differentiation of
$Y_p$ in the trivial bundle $\mathcal{T} \times T_pM_\setC$, and
$\pi_\sigma^{1,0} \colon \mathcal{T} \times TM_\setC \to T_\sigma$ is
the projection.

Now $\hat \nabla^T$ induces a connection $\hat \nabla^K$ in $K =
\bigwedge^m T^*$, which in turn induces a connection $\hat
\nabla^\delta$ in the square root $\delta$. With the help of the
connection $\hat \nabla^\mathcal{L}$, this induces a connection $\hat
\nabla^r$ in the line bundle $\mathcal{\hat L}^k \otimes \delta$.

\begin{definition}
  \label{def:2}
  The connection
  \begin{align*}
    \hat \nabla^r = \big ( \hat \nabla^\mathcal{L} \big )^{\otimes k}
    \otimes \Id + \Id \otimes \hat \nabla^\delta
  \end{align*}
  in $\mathcal{\hat L}^k \otimes \delta \to \mathcal{T} \times M$ is
  called the \emph{reference connection}.
\end{definition}

Notice how the reference connection induces a connection in
$\mathcal{H}^{(k)} \to \mathcal{T}$. Indeed, for any section $s$ of
$\mathcal{H}^{(k)}$ (which is the same as a section of $\mathcal{\hat
  L}^k \otimes \delta$ over ${\mathcal T}\times M$) and any vector field $V$ tangent to
$\mathcal{T}$, it is simply given by $\hat \nabla^r_V s$. Moreover, if
we restrict to a point $\sigma \in \mathcal{T}$ and take $X$ to be a
vector field tangent to $M$, then $(\hat \nabla^r_X s)_\sigma =
(\nabla_\sigma)_X s_\sigma$, so the reference connection is a unified
description of a connection in $\mathcal{H}^{(k)}$ and the connections
in the bundles $\mathcal{L}^k \otimes \delta_\sigma \to M$.

\section*{Curvature}

Later we shall have need for the curvature of the reference
connection, which is given by \Fref{prop:2}, \ref{prop:4} and
\ref{prop:5} below.
\begin{proposition}
  \label{prop:2}
  For vector fields $X$ and $Y$, tangent to $M$, we have
  \begin{align}
    \label{eq:5}
    R_{\hat \nabla^r}(X,Y) = -ik\omega(X,Y) + \frac{i}{2} \rho(X,Y),
  \end{align}
  where $\rho_\sigma$ denotes the Ricci form on $M_\sigma$.
\end{proposition}

\begin{proof}
  This follows immediately by the curvature of prequantum line bundles
  and the standard fact that the canonical line bundle $K_\sigma$ over
  $M_\sigma$ has curvature $i\rho_\sigma$.
\end{proof}

Before giving the curvature in the mixed directions, we introduce some
more notation. Since the symplectic form is non-degenerate, it induces
an isomorphism
\begin{align*}
  i_\omega \colon TM_\setC \to TM_\setC^*,
\end{align*}
by contraction in the first entry.  Moreover $\omega$ is
$J$-invariant, or equivalently of of type (1,1), which implies that
$i_\omega$ interchanges types. Similarly the metric induces a
type-interchanging isomorphism $i_g \colon TM_\setC \to TM_\setC^*$,
and the two are related by $i_{g} = -J i_\omega$.

For any vector field $V$ tangent to $\mathcal{T}$, we can
differentiate the family of complex structures in the direction of $V$
and obtain
\begin{align*}
  V[J] \colon \mathcal{T} \to C^\infty(M, \End(TM_\setC)).
\end{align*}
By differentiation of the identity $J^2 = -\Id$, we see that $V[J]$
anticommutes with $J$. This in turn implies that $V[J]_\sigma$
interchanges types on $M_\sigma$, whence it decomposes as
\begin{align*}
  V[J]_\sigma = V[J]_\sigma' + V[J]_\sigma'',
\end{align*}
where $V[J]_\sigma' \in C^\infty(M, \bar T_\sigma^* \otimes T_\sigma)$
and $V[J]_\sigma'' \in C^\infty(M, T_\sigma^* \otimes \bar T_\sigma)$.

Now define $\tilde G(V) \in C^\infty(M, TM_\setC \otimes TM_\setC)$ by
the relation
\begin{align*}
  V[J] = (\Id \otimes i_\omega)(\tilde G(V))
\end{align*}
for all vector fields $V$. We use the notation
\[\tilde G(V) \omega = (\Id \otimes i_\omega)(\tilde G(V)).\]
The way to interpret this, is to trace the right contravariant part of
$\tilde G(V)$ with the left covariant part of $\omega$, as prescribed
by $(\Id \otimes i_\omega)(\tilde G(V))$. Now observe, that the
combined types of $V[J]$ and $\omega$ yield a decomposition
\begin{align}
  \label{eq:12}
  \tilde G(V) = G(V) + \bar G(V),
\end{align}
for all real vector fields $V$ on $\mathcal{T}$, where $G(V)_\sigma
\in C^\infty(M, T_\sigma \otimes T_\sigma)$ and $\bar G(V)_\sigma \in
C^\infty(M, \bar T_\sigma \otimes \bar T_\sigma)$. Differentiating the
definition, $g = \omega J$, of the metric along $V$ we have
\begin{align}
  \label{eq:13}
  V[g] = \omega V[J] = \omega \tilde G(V) \omega = - (i_\omega \otimes
  i_\omega)( \tilde G(V)).
\end{align}
Once again, notice how the notation $\omega \tilde G(V) \omega$ is
used to denote tracing the right covariant part of $\omega$ with the
left contravariant part of $\tilde G(V)$, as well as tracing the right
contravariant part of $\tilde G(V)$ with the left covariant part of
$\omega$. Since $g$ is symmetric, so is $V[g]$, which implies that
$G(V)_\sigma \in C^\infty(M, S^2(T_\sigma))$ and $\bar G(V)_\sigma \in
C^\infty(M, S^2(\bar T_\sigma))$.\\

By a small calculation, we obtain another useful formula for the
connection $\hat \nabla^T$ in the directions tangent to
$\mathcal{T}$. Indeed, we have that
\begin{equation}
  \label{eq:8}
  \hat \nabla^T_V Y = V [\pi^{1,0} Y] - V[\pi^{1,0}] Y = V[Y] +
  \frac{i}{2} V[J] Y,
\end{equation}
for any section $Y$ of $T$.\\

Now we are ready to calculate the curvature of the reference
connection in the remaining directions. To do this we recall the
general fact, which was already implicitly used to find the curvature
of the half-form bundle, that the curvature of $\hat \nabla^\delta$ is
given by
\begin{align}
  \label{eq:10}
  R_{\hat \nabla^\delta} = -\frac{1}{2} \Tr R_{\hat \nabla^T},
\end{align}
where we trace the endomorphism part of $R_{\hat \nabla^T} \in
\Omega^2(T \times M, \End(T))$. The change of sign appears when we
induce $\hat \nabla^T$ in $T^*$, the trace appears when we induce in
$K = \bigwedge^m T^*$, and the division by two appears when we induce
in $\delta$. Then we have

\begin{proposition}
  \label{prop:4}
  For vector fields $V$ and $W$ tangent to $\mathcal{T}$ we have
  \begin{align}
    \label{eq:11}
    R_{\hat \nabla^r}(V,W) = 0
  \end{align}
\end{proposition}

\begin{proof}
  Take $V$ and $W$ to be pullbacks of vector fields on $\mathcal{T}$
  which satisfy that $[V,W] = 0$. Then using \eqref{eq:8}, we find
  that
  \begin{align*}
    \hat \nabla^T_V \hat \nabla^T_W Y &= \hat \nabla^T_V( W[Y] +
    \frac{i}{2} W[J]Y) \\ &= VW[Y] + \frac{i}{2} VW[J] Y + \frac{i}{2}
    W[J]V[Y] + \frac{i}{2} V[J] W[Y] - \frac{1}{4}V[J]W[J]Y.
  \end{align*}
  Using that $V$ and $W$ commute we get
  \begin{align*}
    R_{\hat \nabla^T}(V,W) Y &= \hat \nabla^T_V \hat \nabla^T_W Y -
    \hat \nabla^T_W \hat
    \nabla^T_V Y \\ &= -\frac{1}{4} (V[J]W[J] - W[J] V[J]) Y \\
    &= -\frac{1}{4}[V[J],W[J]] Y,
  \end{align*}
  and so by \eqref{eq:10} we get
  \begin{align*}
    R_{\hat \nabla^r}(V, W) &= R^{(k)}_{\hat \nabla^\mathcal{L}}(V,W)
    - \frac{1}{2} \Tr R_{\hat \nabla^T}(V,W) = 0,
  \end{align*}
  as desired, since $R_{\hat \nabla^\mathcal{L}}(V,W) = 0$ and
  $R_{\hat \nabla^T}(V,W)$ is a commutator.
\end{proof}

Now we calculate the curvature of the reference connection in the
mixed directions.
\begin{proposition}
  \label{prop:5}
  For vector fields $V$ and $X$, tangent to $\mathcal{T}$ and $M$
  respectively, we have
  \begin{align}
    \label{eq:14}
    R_{\hat \nabla^r}(V,X) = \frac{i}{4} \Tr \tilde \nabla (\tilde
    G(V)) \omega X
  \end{align}
\end{proposition}

\begin{proof}
  First we calculate the curvature of $\hat \nabla^T$. Let $X$ and $V$
  be pullbacks of real vector fields on $M$ and $\mathcal{T}$
  respectively, and let $Y$ be any section of $T$. Then we get
  \begin{align*}
    R_{\hat \nabla^T}(V, X) Y &= \hat \nabla^T_V \hat \nabla^T_X Y -
    \hat \nabla^T_X \hat \nabla^T_V Y \\ &= \pi^{1,0} V[\tilde
    \nabla_X Y] - \tilde \nabla_X \pi^{1,0} V[Y] \\ &= \pi^{1,0}
    V[\tilde \nabla_X Y] - \pi^{1,0} \tilde \nabla_X V[Y] \\ &=
    \pi^{1,0} V[\tilde \nabla]_X Y
  \end{align*}

  By Theorem 1.174 in \cite{MR867684}, we get that the variation of
  the Levi-Civita connection in the tangent bundle is a symmetric
  (2,1)-tensor given by
  \begin{align}
    \label{eq:15}
    \begin{aligned}
      g(V[\tilde \nabla]_X Y, Z) = \tfrac{1}{2} ( &\tilde \nabla_X
      (V[g])(Y,Z) \\
      + &\tilde \nabla_Y (V[g])(X,Z) \\ - &\tilde \nabla_Z(V[g])(X,Y))
    \end{aligned}
  \end{align}
  for vector fields $X$, $Y$ and $Z$ on $M$ and $V$ on $\mathcal{T}$.
  We focus our attention on a point $p \in M$, and let $e_1, \ldots,
  e_m$ be a basis of $T_pM$ satisfying the orthogonality condition
  that $g(e'_j,e''_l)=\delta_{jl}$. Then
  \begin{align*}
    \Tr R_{\hat \nabla^T}(V, X) = \Tr \pi^{1,0} V[\tilde \nabla]_X
    \pi^{1,0} = \sum_\nu g(V[\tilde \nabla]_X e'_\nu, e_\nu'').
  \end{align*}
  But taking into account the type of $V[g]$, and the fact that
  $\tilde \nabla$ preserves types, we get
  \begin{align*}
    g(V[\tilde \nabla]_X e'_\nu, e''_\nu) &= \frac{1}{2} \tilde
    \nabla_{e'_\nu}(V[g])(X,e''_\nu) - \frac{1}{2} \tilde
    \nabla_{e''_\nu} (V[g])(X, e'_\nu) \\ &= \frac{1}{2} X \omega
    \tilde \nabla_{e'_\nu}(\tilde G(V)) \omega e''_\nu - \frac{1}{2} X
    \omega \tilde \nabla_{e''_\nu}(\tilde G(V)) \omega e'_\nu \\ &=
    \frac{i}{2} X \omega \tilde \nabla_{e'_\nu}(G(V)) g e''_\nu +
    \frac{i}{2} X \omega \tilde \nabla_{e''_\nu}(\bar G(V)) g e'_\nu
    \\ &= - \frac{i}{2} g(\tilde \nabla_{e'_\nu}(G(V)) \omega X,
    e''_\nu) - \frac{i}{2} g(\tilde \nabla_{e''_\nu}(\bar G(V)) \omega
    X, e'_\nu).
  \end{align*}
  Summing over $\nu$, we conclude that
  \begin{align*}
    \Tr R_{\hat \nabla^T}(V, X) &= - \frac{i}{2} \Tr \tilde \nabla
    (G(V)) \omega X - \frac{i}{2} \Tr \tilde \nabla (\bar G(V)) \omega
    X \\ &= - \frac{i}{2} \Tr \tilde \nabla(\tilde G(V)) \omega X,
  \end{align*}
  at the point $p$ which was arbitrary.  Finally we get by
  \Fref{lem:1} and \eqref{eq:10} that
  \begin{align*}
    R_{\hat \nabla^r}(V, X) &= R^{(k)}_{\hat \nabla^\mathcal{L}}(V,X)
    - \frac{1}{2} \Tr R_{\hat \nabla^T}(V,X) \\ &= \frac{i}{4} \Tr
    \tilde \nabla(\tilde G(V)) \omega X,
  \end{align*}
  which was the claim.
\end{proof}

\chapter{The Hitchin Connection}

Let $\mathcal{D}(M_\sigma, \mathcal{L}^k \otimes \delta_\sigma)$
denote the space of differential operators on
$\mathcal{H}^{(k)}_\sigma = C^\infty(M_\sigma, \mathcal{L}^k \otimes
\delta_\sigma)$, and consider the bundle $\mathcal{D}(M, \mathcal{\hat
  L}^k \otimes \delta)$ over $\mathcal{T}$ having these spaces as
fibers. One could think of $\mathcal{D}(M, \mathcal{\hat L}^k \otimes
\delta)$ as the space of differential operators on sections of
$\mathcal{\hat L}^k \otimes \delta$, which are of order zero in the
directions tangent to $\mathcal{T}$. Then, for any one-form $u$ on
$\mathcal{T}$ with values in $\mathcal{D}(M, \mathcal{\hat L}^k
\otimes \delta)$, we have a connection $\boldnabla$ in the bundle
$\mathcal{H}^{(k)} = C^\infty(M, \mathcal{L}^k \otimes \delta)$ over
$\mathcal{T}$ given by
\begin{align*}
  \boldnabla_V = \hat \nabla^r_V + u(V),
\end{align*}
for any vector field $V$ on $\mathcal{T}$. Now we wish to find a $u$
such that $\boldnabla$ preserves the subspaces $H_{\sigma}^{(k)}$,
thereby proving that these form a subbundle and inducing a connection
in this subbundle.

\begin{lemma}
  \label{lem:3}
  The connection $\boldnabla$ preserves $H^{(k)}$ if and only if
  \begin{align}
    \label{eq:16}
    \nabla^{0,1} u(V) s = \frac{i}{2} V[J] \nabla s + \frac{i}{4} \Tr
    \tilde \nabla (G(V)) \omega s,
  \end{align}
  for all vector fields $V$ on $\mathcal{T}$, and all $s \in H^{(k)}$.
\end{lemma}

\begin{proof}
  Let $X$ and $V$ be the pullbacks of a vector field on $M$ and
  $\mathcal{T}$ respectively. Then we see that
  \begin{align}
    \label{eq:17}
    [V, X''] = \frac{i}{2} V[J] X.
  \end{align}
  Now, assume that $s \in H^{(k)}_\sigma$ and consider any extension
  of $s$ to a smooth section of $\mathcal{H}^{(k)} \to
  \mathcal{T}$. Then we get
  \begin{align*}
    \nabla_{X''} \boldnabla_V s &= \hat \nabla^r_{X''} \hat \nabla^r_V
    s + \nabla_{X''} u(V)s \\ &= \hat \nabla^r_V \hat \nabla^r_{X''} s
    - R_{\hat \nabla^r}(V,X'') s - \hat \nabla^r_{[V,X'']} s +
    \nabla_{X''} u(V)s \\
    &= - \frac{i}{2} \nabla_{V[J] X} s - \frac{i}{4} \Tr (\tilde
    \nabla (G(V)) \omega X)s + \nabla_{X''} u(V)s,
  \end{align*}
  at the point $\sigma \in \mathcal{T}$, where we used \eqref{eq:17}
  and \Fref{prop:5} for the last equality.  This tells us, that
  $\boldnabla$ preserves $H^{(k)}$ if and only if $u$ satisfies the
  equation in the lemma.
\end{proof}

For any vector field $V$ tangent to $\mathcal{T}$, the tensor
$G(V)_\sigma \in C^\infty(M_\sigma,S^2(T_\sigma))$ induces a linear
map $G(V)_\sigma \colon TM_\setC^* \to TM_\setC$, by the formula
\begin{align*}
  \alpha \mapsto \Tr(G(V)_\sigma \otimes \alpha) = G(V)_\sigma\alpha.
\end{align*}
Obviously this is in fact a map $G(V)_\sigma \colon T_\sigma^*\to
T_\sigma$. We then define a second order operator
$\Delta_{G(V)_\sigma} \in \mathcal{D}(M, \mathcal{L}^k \otimes
\delta_\sigma)$ by $\Delta_{G(V)_\sigma} = \Tr \nabla_\sigma
G(V)_\sigma \nabla_\sigma$, or more explicitly by the diagram
\begin{align}
  \label{eq:18}
  \begin{aligned}
    \xymatrix@C=50pt@R=20pt{ & {C^{\infty}(M_\sigma,TM_\setC^* \otimes
        \mathcal{L}^k\otimes \delta_\sigma)}
      \ar[d]^-{G(V)_\sigma\otimes\id\otimes\id} \\
      {C^{\infty}(M_\sigma, \mathcal{L}^k\otimes \delta_\sigma)}
      \ar[]+UR;[ur]+L^-{\nabla_\sigma}&
      {C^{\infty}(M_\sigma,T_\sigma\otimes \mathcal{L}^k\otimes
        \delta_\sigma)} \ar[d]^-{\tilde \nabla_\sigma \otimes \id +
        \id\otimes
        \nabla_\sigma}\\
      & {C^{\infty}(M_\sigma,TM_\setC^* \otimes T_\sigma\otimes
        \mathcal{L}^k\otimes \delta_\sigma)} \ar[]+L;[ul]+DR_-{\Tr} }
  \end{aligned}
\end{align}

We shall make the additional assumption, that the family $J$ is rigid
in the sense that $G(V)_\sigma$ should be a holomorphic section of
$S^2(T_\sigma)$ over $M_\sigma$.
\begin{definition}
  \label{def:4}
  The family $J$ of k\"ahler structures on $(M,\omega)$ is called
  \emph{rigid} if
  \begin{align*}
    \tilde \nabla^{0,1}_\sigma(G(V)_\sigma) = 0
  \end{align*}
  for all vector fields $V$ tangent to $\mathcal{T}$ and $\sigma \in
  \mathcal{T}$.
\end{definition}
From now on, we will for simplicity often suppress the subscription $\sigma$
from the notation.  Under this assumption we have the following lemma
\begin{lemma}
  \label{lem:4}
  At every point $\sigma \in \mathcal{T}$, the operator $\Delta_{G(V)}
  = \Tr \nabla G(V) \nabla$ satisfies
  \begin{align*}
    \nabla^{0,1}\Delta_{G(V)} s = - 2ik \omega G(V) \nabla s + ik\Tr
    \tilde \nabla (G(V)) \omega s - \tfrac{i}{2} \Tr \tilde \nabla
    (G(V) \rho) s
  \end{align*}
  for all vector fields $V$ tangent to $\mathcal{T}$ and all (local)
  holomorphic sections $s$ of the line bundle $\mathcal{L}^k \otimes
  \delta \to M$.
\end{lemma}

\begin{proof}
  The proof is by direct calculation. Letting $G$ denote $G(V)$ we
  have
  \begin{align*}
    \nabla^{0,1} \Delta_{G} s = \nabla^{0,1} \Tr \nabla G \nabla s =
    \Tr\nabla^{0,1} \nabla G \nabla s.
  \end{align*}
  Working further on the right side we commute the two connections,
  giving as ekstra terms the curvature of $M_\sigma$ and of the line
  bundle $\mathcal{L}^k \otimes \delta_\sigma$,
  \begin{align*}
    \nabla^{0,1} \Delta_{G} s = \Tr\nabla \nabla^{0,1} G \nabla s - ik
    \omega G \nabla^{1,0} s + \tfrac{i}{2} \rho G \nabla s - i \rho G
    \nabla s.
  \end{align*}
  Collecting the last two terms, and using the fact that $J$ is rigid
  on the first, we obtain
  \begin{align*}
    \nabla^{0,1} \Delta_{G} s = \Tr\nabla G \nabla^{0,1} \nabla s - ik
    \omega G \nabla s - \tfrac{i}{2} \rho G \nabla s.
  \end{align*}
  Commuting the two connections, and using that $s$ is holomorphic, we
  get
  \begin{align*}
    \nabla^{0,1} \Delta_{G} s = ik \Tr \nabla G \omega s -
    \tfrac{i}{2} \Tr \nabla G \rho s - ik \omega G \nabla s -
    \tfrac{i}{2} \rho G \nabla s.
  \end{align*}
  Expanding the covariant derivatives in the first two terms by the
  Leibniz rule, and using the fact that $\omega$ is parallel, we get
  the following, after collecting and cancelling terms
  \begin{align*}
    \nabla^{0,1} \Delta_{G} s = ik \Tr \tilde \nabla (G) \omega s -
    2ik \omega G \nabla s - \tfrac{i}{2}\Tr \tilde \nabla (G \rho)s
  \end{align*}
  This was the desired expression. Moreover we notice, that the above
  is a local computation, so that the identity is valid for local
  holomorphic sections of $\mathcal{L}^k \otimes \delta$ as well.
\end{proof}

\begin{corollary}
  \label{cor:1}
  Provided that $H^{0,1}(M) = 0$, we have that $\Tr \tilde \nabla
  (G(V) \rho)$ is exact with respect to the $\bar \d$-operator on
  $M$.
\end{corollary}
\begin{proof}
  By appealing to \Fref{lem:4}, in the case where $k = 0$, we get for
  any local holomorphic section $s$ of $\mathcal{L}^k \otimes
  \delta_\sigma \to M_\sigma$ that
  \begin{align*}
    0 = \frac{i}{2} \nabla^{0,1}_\sigma \Tr \tilde \nabla_{\sigma}
    (G(V)_\sigma \rho_\sigma)s = \frac{i}{2} \bar \d_\sigma ( \Tr
    \tilde \nabla_{\sigma} (G(V)_\sigma \rho_\sigma)) s.
  \end{align*}
  This
  immediately implies that
  \begin{align*}
    0 = \bar \d_\sigma ( \Tr \tilde \nabla_{\sigma} (G(V)_\sigma
    \rho_\sigma)),
  \end{align*}
  and since $H^{0,1}(M) = 0$, the corollary follows.
\end{proof}

We remark, that the assumption $H^{0,1}(M)=0$ is satisfied for any
compact K\"ahler manifold with $H^1(M, \setR) = 0$, by the Hodge
decomposition theorem.

By \Fref{cor:1}, we choose any smooth one-form $H \in
\Omega^1(\mathcal{T}, C^\infty(M))$, such that
\begin{align}
  \label{eq:19}
  \bar \d H(V) = \frac{i}{2} \Tr \tilde \nabla (G(V) \rho),
\end{align}
for any vector field $V$ on $\mathcal{T}$. Then finally we define
\begin{equation}
  \label{eq:20}
  u(V) = \frac{1}{4k}(\Delta_{G(V)} + H(V)),
\end{equation}
which clearly solves equation \eqref{eq:16}. Thus we have proved
\Fref{thm:1}.

\chapter{Relation to Non-Corrected Quantization}

We now impose the same assumptions as in \cite{A1} in order to
give an explicit solution and to compare the constructed Hitchin
connection with that previously constructed in \cite{A1}.

Thus, from now on $M$ is assumed to be compact with $H^1(M, \setR) =
0$. The real first Chern class of $(M, \omega)$, that is the image of
the first Chern class in $H^2(M,\setR)$, is assumed to satisfy
\begin{align}
  \label{eq:21}
  c_1(M, \omega) = n[\tfrac{\omega}{2\pi}],
\end{align}
where $n \in \setZ$ is some integer, which must be even by our
assumption on the second Stiefel-Whitney class of $M$. Finally
$\mathcal{T}$ is assumed to be a complex manifold and the map $J$ to
be \emph{holomorphic} in the following sense.

\begin{definition}
  \label{def:3}
  The family $J$, of K\"ahler structures on $M$ parametrized by
  $\mathcal{T}$, is called \emph{holomorphic} if it satisfies
  \begin{align*}
    V'[J] = V[J]' \qquad \text{and} \qquad V''[J] = V[J]'',
  \end{align*}
  for every vector field $V$ tangent to $\mathcal{T}$.
\end{definition}

These assumptions have a number of consequences which we shall explore
in the following. First, we give an alternative characterization of
holomorphic families of K\"ahler structures.

Let $I$ denote the integrable almost complex structure on
$\mathcal{T}$ induced by its complex structure. Then we have an almost
complex structure $\hat J$ on $\mathcal{T} \times M$ defined by
\begin{align}
  \label{eq:6}
  \hat J(V \oplus X) = IV \oplus J_\sigma X, \qquad V \oplus X \in
  T_{(\sigma,p)} (\mathcal{T} \times M).
\end{align}
The following gives another characterization of holomorphic
families.
\begin{proposition}
  \label{prop:3}
  The family $J$ is holomorphic if and only if $\hat J$ is integrable.
\end{proposition}

\begin{proof}
  We show that $J$ is holomorphic if and only if the Nijenhuis tensor
  for $\hat J$ vanishes. By the Newlander-Nirenberg theorem this will
  imply the proposition (See e.g. \cite{MR1393941}).

  Clearly the Nijenhuis tensor vanishes, when evaluated only on
  vectors tangent to $\mathcal{T}$, since $I$ is integrable. Likewise
  it will vanish when evaluated only on vectors tangent to $M$, since
  $J$ is a family of integrable almost complex structures. Thus we are
  left with the case of mixed directions.

  Let $X$ and $V$ be pullbacks to $\mathcal{T} \times M$ of vector
  fields on $M$ and $\mathcal{T}$ respectively. Then since $X$ is
  constant along $\mathcal{T}$ and $V$ is constant along $M$ we find
  that
  \begin{align}
    \label{eq:7}
    [V, JX] = V[J] X.
  \end{align}
  Now consider the following evaluation of the Nijenhuis tensor
  \begin{align*}
    N(V', X) &= [IV', JX] - [V', X] - \hat J[IV', X] - \hat J[V', JX]
    \\ &= i[V', JX] - \hat J[V', JX] \\ &= iV'[J] X - J V'[J]X \\ &=
    2i \pi^{0,1} V'[J] X.
  \end{align*}
  Similarly one shows, that $N(V'', X) = -2i \pi^{1,0} V''[J] X$. Thus
  we see that $N(V, X)$ vanishes if and only if
  \begin{align*}
    \pi^{0,1} V'[J] X = 0 \quad \text{and} \quad \pi^{1,0} V''[J] X =
    0.
  \end{align*}
  This proves the proposition.
\end{proof}

We shall denote by $\hat d$ the differential on $\mathcal{T} \times
M$, which splits as
\begin{align*}
  \hat d = d_\mathcal{T} + d_M
\end{align*}
into the sum of the differentials on $\mathcal{T}$ and $M$
respectively. Similar notation is used for $\partial$ and $\bar\partial$.

\subsection*{Explicit Formula for $H(V)$}
\label{sec:explicit-formula}

As a first consequence of our additional assumptions we are able to give
an explicit formula for the one-form $H$ in \eqref{eq:20}.

Since the curvature of the canonical line bundle $K_\sigma$ is
$i\rho_\sigma$, the real first Chern class of $M_\sigma$ is
represented by $\tfrac{\rho_\sigma}{2\pi}$. Since the K\"ahler form is
harmonic, the assumption \eqref{eq:21} is then equivalent to
$\rho^H_\sigma = n\omega$, where $\rho^H_\sigma$ denotes the harmonic
part of the Ricci form.

Since any real exact (1,1)-form on a K\"ahler manifold is $\d \bar
\d$-exact, there exists, for any $\sigma \in \mathcal{T}$, a real
function $F_\sigma$, called a \emph{Ricci potential}, satisfying
\begin{align*}
  \rho_\sigma = \rho^H_\sigma + 2i\d_\sigma \bar \d_\sigma F_\sigma.
\end{align*}
By compactness of $M$, any two Ricci potentials on $M_\sigma$ differ
by a constant. Thus choosing a particular normalization, such as
\begin{align}
  \label{eq:22}
  \int_M F_\sigma \omega^m = 0,
\end{align}
would yield a real smooth function $F \in C^\infty(\mathcal{T} \times
M)$, with $F_\sigma$ a Ricci potential on $M_\sigma$ for every $\sigma
\in \mathcal{T}$. Such a function shall be called a \emph{smooth family of
Ricci potentials over $\mathcal{T}$}, and it satisfies the identity
\begin{align}
  \label{eq:23}
  \rho = n\omega + 2i\d_M \bar \d_M F.
\end{align}
By a \emph{pluriharmonic family of Ricci potentials over
  $\mathcal{T}$} we mean a smooth family of Ricci potentials satisfying
\begin{align}
  \label{eq:9}
  \d_\mathcal{T} \bar \d_\mathcal{T} F = 0.
\end{align}
We shall have use for this notion later, but for the moment we just
consider $F$ to be any smooth family of Ricci potentials.

We will need the following lemma, the proof of which is given in
\cite{A1}.
\begin{lemma}
  \label{lem:5}
  Any smooth family $F$ of Ricci potentials satisfies
  \begin{align}
    \label{eq:24}
    \bar \d_M V'[F] = - \frac{i}{4} \Tr \tilde \nabla (G(V)) \omega -
    \frac{i}{2} \d_M F G(V) \omega,
  \end{align}
  for any vector field $V$ tangent to $\mathcal{T}$.
\end{lemma}

Then we have the following

\begin{lemma}
  \label{lem:6}
  Let $F$ be a smooth family of Ricci potentials. Then the
  one-form $H \in \Omega^1(\mathcal{T}, C^\infty(M))$ given by
  \begin{align*}
    H(V) = - 2nV'[F] - \d_M F G(V) \d_M F - \Tr \tilde \nabla ( G(V)
    \d_M F)
  \end{align*}
  satisfies $\bar \d_M H(V) = \frac{i}{2} \Tr \tilde \nabla
  (G(V)\rho)$.
\end{lemma}

\begin{proof}
  Throughout this proof we shall denote $\d_M$ and $\bar \d_M$ for
  short by $\d$ and $\bar \d$ respectively. Since $\omega$ is
  parallel, with respect to the Levi-Civita connection $\tilde
  \nabla$, we get
  \begin{align*}
    \Tr \tilde \nabla (G(V) \rho) &= \Tr
    \tilde \nabla (G(V) (n\omega + 2i\d \bar \d F)) \\ &=
    n\Tr\tilde \nabla (G(V)) \omega + 2i \Tr \tilde \nabla (G(V) \d \bar
    \d F).
  \end{align*}
  Moreover, it is easily verified that
  \begin{align*}
    \Tr \tilde \nabla (G(V) \d \bar \d F) &= - i \d F G(V) \rho + \bar
    \d \Tr \tilde \nabla (G(V) \d F) \\ &= - in \d F G(V) \omega + 2 \d
    F G(V) \d \bar \d F + \bar \d \Tr \tilde \nabla ( G(V) \d F ).
  \end{align*}
  Then the lemma follows by \Fref{lem:5} and the identity
  \begin{align*}
    \bar \d (\d F G(V) \d F) = 2 \d F G(V) \d \bar \d F,
  \end{align*}
  which is easily verified, using the symmetry of $G(V)$.
\end{proof}

Thus under the assumptions of this section, we have a completely explicit
formula for the Hitchin connection.\\

\subsection*{Curvature of the Reference Connection Revisited}
\label{sec:curvature-revisited}

Notice that the type of $\omega$, and the fact that $J$ is
holomorphic, implies
\begin{align*}
  V'[J] = V[J]' = G(V) \omega,
\end{align*}
which in turn gives $G(V) = G(V')$.  Then, having calculated the
curvature of the reference connection in all directions, we see that
it is of type (1,1) over $\mathcal{T} \times M$ and thus the
(0,2)-part of the curvature vanishes.  This means that the reference
connection defines a holomophic structure on the line bundle
$\mathcal{\hat L}^k \otimes \delta$, over the complex manifold
$\mathcal{T} \times M$.  Moreover, we observe that $(\hat
\nabla^r)^{0,1}$ preserves the bundle $H^{(k)} \to \mathcal{T}$, since
$u(V'') = 0$ solves \eqref{eq:16}. Thus the reference connection
defines a holomorphic
structure on the bundle $H^{(k)} \to \mathcal{T}$.\\

We now prove that, at least locally over $\mathcal{T}$, the curvature
of the reference connection can be expressed in terms of a
pluriharmonic family of Ricci potentials.

First we have the following lemma, which is an immediate consequence
of \Fref{lem:5}, by direct verification.
\begin{proposition}
  \label{prop:6}
  For any smooth family $F$ of Ricci potentials and any vector fields
  $V$ on $\mathcal{T}$ and $X$ on $M$, the curvature of the reference
  connection is given by
  \begin{align*}
    R_{\hat \nabla^r}(V, X) = -\hat \d \bar{ \hat \d} F (V, X)
  \end{align*}
\end{proposition}

\begin{proof}
  Let $V$ and $X$ be pullbacks of real vector fields on $\mathcal{T}$
  and $M$ respectively. Then we have
  \begin{align*}
    \bar{\hat \d} \hat \d F(X'', V') &= \hat d \hat \d F (X'', V') \\
    &= X'' (\hat \d F(V')) - V'( \hat \d F(X'')) - \hat \d F
    ([X'', V']) \\ &= X'' V'[F] + \frac{i}{2} \hat \d F (V'[J] X) \\
    &= X'' V'[F] + \frac{i}{2} \d_M F G(V) \omega X'' \\ &=
    -\frac{i}{4} \Tr \tilde \nabla (G(V)) \omega X'' \\ &= -R_{\hat
      \nabla^r}(V', X''),
  \end{align*}
  where we use \Fref{lem:5} and \Fref{prop:5} for the last two
  equalities. The case of $X'$ and $V''$ is similar by conjugation of
  the identity in \Fref{lem:5}.
\end{proof}

Then we have

\begin{theorem}
  \label{thm:3}
  Let $(M, \omega)$ be a compact, prequantizable, symplectic manifold
  with the real first Chern class satisfying $c_1(M, \omega) =
  n[\tfrac{\omega}{2\pi}]$, $H^1(M, \setR) = 0$ and vanishing second
  Stiefel-Whitney class. Let $J$ be a rigid, holomorphic family of
  K\"ahler structures on $M$, parametrized by a complex manifold
  $\mathcal{T}$. Then for every pluriharmonic family of Ricci
  potentials $\tilde F$ we have
  \begin{align}
    \label{eq:36}
    R^{(k)}_{\hat \nabla^r} = R^{(k-n/2)}_{\hat \nabla^\mathcal{L}} -
    \hat \d \bar{\hat \d} \tilde F,
  \end{align}
  where $R^{(k)}_{\hat \nabla^r}$ denotes curvature of the reference
  connection in $\mathcal{\hat L}^k \otimes \delta$ and
  $R^{(k-n/2)}_{\hat \nabla^\mathcal{L}}$ denotes the curvature of
  $\hat \nabla^\mathcal{L}$ in $\mathcal{\hat L}^{k - n/2}$.
\end{theorem}
\begin{proof}
  Let $X$ and $Y$ be vector fields tangent to $M$, and let $V$ and $W$
  be vectorfields tangent to $\mathcal{T}$. Then by \Fref{prop:2} and
  \eqref{eq:23} we have that
  \begin{align*}
    R_{\hat \nabla^r}(X,Y) &= -ik\omega(X,Y) + \frac{i}{2} \rho(X,Y) \\
    &= -i(k-\tfrac{n}{2}) \omega(X,Y) - \d_M \bar \d_M \tilde F (X,Y)
    \\ &= R^{(k-n/2)}_{\hat \nabla^\mathcal{L}} (X,Y) - \hat \d \bar{
      \hat \d} \tilde F (X,Y).
  \end{align*}
  By \Fref{lem:1}, the curvature $R^{(k-n/2)}_{\hat
    \nabla^\mathcal{L}}$ vanishes in the remaining directions, and so
  the theorem follows from from \Fref{prop:6} and the fact that
  $\tilde F$ is a pluriharmonic family.
\end{proof}

There is no guarantee that pluriharmonic families of Ricci potentials
exist globally over $\mathcal{T}$. However we are able to prove that
at least locally they do exist. First we prove

\begin{lemma}
  \label{lem:7}
  For any smooth family $F$ of Ricci potentials, and any vector fields
  $V$ and $W$ on $\mathcal{T}$, we have
  \begin{align*}
    0 = d_M \big [\hat \d \bar{\hat \d} F
    (V, W) \big ].
  \end{align*}
\end{lemma}

\begin{proof}
  The first equality is just \Fref{prop:4} so we shall prove the
  second. Take $V$, $W$ and $X$ to be commuting vector fields so that
  $V$ and $W$ are tangent to $\mathcal{T}$ and $X$ is tangent to
  $M$. Then we must prove
  \begin{align*}
    0 = X\big [ \hat \d \bar {\hat \d} F(V, W) \big ].
  \end{align*}
  Now by the differential Bianchi identity and \Fref{prop:4} we have
  \begin{align*}
    0 &= d^{\nabla^r} R_{\nabla^r}(X, W, V) \\ &=
    \vphantom{\frac{1}{1}} \nabla^r_{X} R_{\nabla^r}(W, V) -
    \nabla^r_{W} R_{\nabla^r}(X, V) + \nabla^r_{V} R_{\nabla^r}(X, W)
    \\ &= \nabla^r_{V} R_{\nabla^r}(X, W) - \nabla^r_{W}
    R_{\nabla^r}(X, V).
  \end{align*}
  Then \Fref{prop:6} yields
  \begin{align*}
    0 &= V[\hat \d \bar {\hat \d}F(X,W)] - W[\hat \d \bar {\hat
      \d}F(X,V)] \\ &= VXW''[F] - VWX''[F] - WXV''[F] + WVX''[F] \\ &=
    XVW''[F] - XWV''[F] \\ &= X[\hat \d \bar {\hat \d}F(V, W)]
  \end{align*}
  as desired.
\end{proof}

This allows us to prove

\begin{proposition}
  \label{prop:7}
  Around every point $\sigma \in \mathcal{T}$ there is an
  neighbourhood $U$ and a pluriharmonic family $\tilde F$ of Ricci
  potentials over $U$.
\end{proposition}

\begin{proof}
  Let $\sigma \in \mathcal{T}$ and fix a smooth family $F$ of Ricci
  potentials, say the one satisfying \eqref{eq:22}. Let $V$ and $W$ be
  vectorfields tangent to $\mathcal{T}$. Then by \Fref{lem:7} if
  follows, that there exists a 2-form $\alpha \in
  \Omega^2(\mathcal{T})$ on $\mathcal{T}$ such that
  \begin{align}
    \label{eq:35}
    0 = \hat \d \bar{\hat \d} F(V,W) +
    \pi^*_\mathcal{T} \alpha(V,W),
  \end{align}
  where $\pi_\mathcal{T} \colon \mathcal{T} \times M \to \mathcal{T}$
  is the projection. Now from \eqref{eq:35} it is easily shown that
  $\alpha \in \Omega^2(\mathcal{T})$ is closed.  Since $\alpha$ is
  also of type (1,1), there exists a neighbourhood $U$ around $\sigma$
  and a real function $A \in C^\infty(U, \setR)$ such that $\alpha
  \vert_U = \d_\mathcal{T} \bar \d_\mathcal{T} A$. But then we can
  define another smooth family of Ricci potentials over $U$ by
  \begin{align*}
    \tilde F = F \vert _U + A,
  \end{align*}
  and we see that
  \begin{align*}
    0= \hat \d \bar{\hat \d} \tilde F(V,W),
  \end{align*}
  as desired.
\end{proof}

Thus, locally over $\mathcal{T}$, we can express the curvature of the
reference connection by \eqref{eq:36}. Using this result, we are able
to relate our construction of the Hitchin connection to a construction
of Andersen (\cite{A1}) in the non-corrected setting.

\subsection*{Hitchin's Connection in Non-Corrected Quantization}

We wish to relate the quantum spaces of half-form corrected
quantization to the spaces of non-corrected geometric quantization,
with the intent to describe, in the non-corrected setting, our
construction of a Hitchin connection and relate it to the construction
in \cite{A1}.

It turns out, that the choice of prequantum line bundle plays a role
in this. This is because of the choice of metaplectic structure we
made. We note that what we really chose was a half of $c_1(M,\omega)$,
so all we know is that $\delta$ is a line bundle satisfying
$2c_1(\delta)= - c_1(M,\omega)$. So if we impose on $(M,\omega)$ that $n$
divides $c_1(M,\omega)$, we get that $\frac{n}{2}$ divides
$c_1(\delta)$. We will need that the prequantum line bundle is related
to the metaplectic structure in a certain way, and the following lemma
ensures that this is possible.

\begin{lemma}
  \label{lem:2}
  If $c_1(M, \omega)$ is divisible by $n$ in $H^2(M, \setZ)$, there
  exists a prequantum line bundle $\mathcal{L}$ over $M$ such that
  \begin{align*}
    \tfrac{n}{2} c_1(\mathcal{L}) = - c_1(\delta).
  \end{align*}
\end{lemma}
\begin{proof}
  Let $\mathcal{L}_0$ be any prequantum line bundle on $M$ and pick an
  auxiliary K\"ahler structure $J$ on $M$. Let $F_J$ be a Ricci
  potential on $M$ and consider the line bundles
  $(\mathcal{L}_0^{-n/2}, e^{F_J} h^{\mathcal{L}_0})$ and $(\delta_J,
  h^\delta_J)$ over $M$. Then it is easily calculated, that the line
  bundles have the same curvature. Thus, the tensor product of the
  former with the dual of the latter yields a flat Hermitian line
  bundle $L_1$. Since $c_1(\delta)$ is divisible by
  $\tfrac{n}{2}$, there exists a flat Hermitian line bundle
  $L_2$ such that $L_2^{n/2} \cong
  L_1$. Finally the line bundle $\mathcal{L} = \mathcal{L}_0
  \otimes L_2$ has the structure of a prequantum line
  bundle, and $\tfrac{n}{2}c_1(\mathcal{L}) = c_1(\mathcal{L}^{n/2}) =
  - c_1(\delta)$. Thus $\mathcal{L}$ is the desired prequantum line
  bundle.
\end{proof}
From now on, we will assume that our prequantum line bundle satisfies
$\tfrac{n}{2} c_1(\mathcal{L}) = - c_1(\delta)$.  We note, that only
when $H^2(M,\setZ)$ has torsion, is the assumption a further
restriction on $(M,\omega)$, as otherwise the curvature
determines the line bundle completely.\\

Next, let $\tilde F$ be a local pluriharmonic family of Ricci
potentials over $U$, with $H^1(U) = 0$, such that \eqref{eq:36} is
satisfied. We wish to construct an isomorphism $\hat \phi$ of
holomorphic Hermitian line bundles over $U \times M$
\begin{align}
  \label{eq:44}
  \hat \phi \colon (\mathcal{\hat L}^{k-n/2}, e^{\tilde F} \hat
  h^\mathcal{L}) \to (\mathcal{\hat L}^k \otimes \delta, \hat h).
\end{align}
Since $\frac{n}{2}c_1(\mathcal{L}) = -c_1(\delta)$, the line bundles are
isomorphic as complex line bundles. The obstruction to finding the
structure preserving isomorphism $\hat \phi$ lies in the first
cohomology of $U \times M$. But this is trivial by the K\"unneth
formula, since $H^1(U) = 0$ and $H^1(M) = 0$ by assumption.

Moreover, it is easily seen that the pullback under $\hat \phi$ of the
reference connection is given by
\begin{align}
  \label{eq:37}
  \hat \phi^* \hat \nabla^r = \hat \nabla^\mathcal{L} + \hat \d \tilde
  F,
\end{align}
since the right hand side is the unique Hermitian connection
compatible with the holomorphic structure of $\mathcal{\hat L}^{k-n/2}$.\\

In the paper \cite{A1}, Andersen constructs a Hitchin connection in
the bundle $\mathcal{T} \times C^\infty(M, \mathcal{L}^k)$, preserving
the subbundle of holomorphic sections. His construction is valid for
any rigid, holomorphic family of K\"ahler structures on $M$
parametrized by $\mathcal{T}$, provided that $H^1(M, \setR) = 0$ and
$c_1(M, \omega) = n[\tfrac{\omega}{2\pi}]$.

Now, the existence of the isomorphism \eqref{eq:44} enables us to
compare his construction to the one presented in this paper. Thus we
shall briefly recall that the Hitchin connection constructed in
\cite{A1} is given by
\begin{align}
  \label{eq:38}
  \tilde \boldnabla_V = \hat \nabla^\mathcal{L}_V + \tilde u(V),
\end{align}
where
\begin{align}
  \label{eq:39}
  \tilde u(V) = \frac{1}{4k + 2n}(\Delta^\mathcal{L}_{G(V)} + 2
  \nabla^\mathcal{L}_{G(V) \d \tilde F} + 4k V'[\tilde F]),
\end{align}
and $\Delta^\mathcal{L}_{G(V)}$ is the operator given by the diagram
\begin{align}
  \label{eq:40}
  \begin{aligned}
    \xymatrix@C=50pt@R=20pt{ & {C^{\infty}(M_\sigma,TM_\setC^* \otimes
        \mathcal{L}^k)}
      \ar[d]^-{G(V)_\sigma \otimes \id} \\
      {C^{\infty}(M_\sigma, \mathcal{L}^k)}
      \ar[]+UR;[ur]+L^-{\nabla^\mathcal{L}}& {C^{\infty}(M_\sigma,
        T_\sigma \otimes \mathcal{L}^k)} \ar[d]^-{\tilde \nabla_\sigma
        \otimes \id + \id\otimes
        \nabla^\mathcal{L}}\\
      & {C^{\infty}(M_\sigma,TM_\setC^* \otimes T_\sigma \otimes
        \mathcal{L}^k)} \ar[]+L;[ul]+DR_-{\Tr} }
  \end{aligned}
\end{align}
We leave it to the reader to verify, using \eqref{eq:37}, that the
pullback by $\hat \phi$ of the operator $\Delta_{G(V)}$, acting on
sections of $\mathcal{\hat L}^k \otimes \delta$, is given by
\begin{align}
  \label{eq:41}
  \hat \phi^* \Delta_{G(V)} = \Delta^\mathcal{L}_{G(V)} + 2
  \nabla^\mathcal{L}_{G(V) \d_M \tilde F} - H(V) - 2nV'[\tilde F],
\end{align}
where $H(V)$ is given by the expression in \Fref{lem:6}, but in terms
of $\tilde F$.

Furthermore, in the bundle $\mathcal{\hat L}^{k-n/2}$, the formula
\eqref{eq:39} becomes
\begin{align}
  \label{eq:42}
  \begin{aligned}
    \tilde u(V) &= \frac{1}{4k}(\Delta^\mathcal{L}_{G(V)} + 2
    \nabla^\mathcal{L}_{G(V) \d_M \tilde F} - 2nV'[\tilde F]) +
    V'[\tilde F] \\
    &= \frac{1}{4k}(\hat \phi^* \Delta_{G(V)} + H(V)) + V'[\tilde F] \\
    &= \hat \phi^* u(V) + V'[\tilde F].
  \end{aligned}
\end{align}
But this means, that the pullback of our Hitchin connection by $\hat
\phi$ is given by
\begin{align}
  \label{eq:43}
  \begin{aligned}
    \hat \phi^* \boldnabla_V &= \hat \phi^* \hat \nabla^r_V + \hat
    \phi^* u(V) \\ &= \hat \nabla^\mathcal{L}_V + V'[\tilde F] + \hat
    \phi^* u(V) \\ &= \hat \nabla^\mathcal{L}_V + \tilde u(V) \\ &=
    \tilde \boldnabla_V
  \end{aligned}
\end{align}
Thus the two connections agree, and we have proved \Fref{thm:2}. \\

\clearpage

\renewcommand{\bibname}{References}





\def\cprime{$'$} \def\cprime{$'$} \def\cprime{$'$}

\end{document}